\documentclass{scrartcl}
\usepackage{mathptmx}
\usepackage[latin1]{inputenc}
\usepackage{amsmath}
\usepackage{amsfonts}
\usepackage{amssymb}
\usepackage{subfig}
\usepackage{booktabs}
\usepackage{hyperref}
\usepackage{tikz}
\usetikzlibrary{calc,arrows,decorations.markings,patterns,external}
\newcommand{\lapack}{\textsf{LAPACK}}
\allowdisplaybreaks
\newcommand{\bbbn}{{\mathbb N}}
\newcommand{\bbbr}{{\mathbb R}}
\newcommand{\Idx}{{\mathcal I}}
\newcommand{\Jdx}{{\mathcal J}}
\newcommand{\ctI}{{\mathcal T}_{\Idx}}
\newcommand{\ctJ}{{\mathcal T}_{\Jdx}}
\newcommand{\lfI}{{\mathcal L}_{\Idx}}
\newcommand{\lfJ}{{\mathcal L}_{\Jdx}}
\newcommand{\ctIJ}{{\mathcal T}_{\Idx\times\Jdx}}
\newcommand{\lfIJ}{{\mathcal L}_{\Idx\times\Jdx}}
\newcommand{\lfaIJ}{{\mathcal L}^+_{\Idx\times\Jdx}}
\newcommand{\lfiIJ}{{\mathcal L}^-_{\Idx\times\Jdx}}

\newcommand{\ctt}{{\mathcal T}_{\hat t_0}}
\newcommand{\cts}{{\mathcal T}_{\hat s_0}}

\newcommand{\nI}{n_{\Idx}}
\newcommand{\nJ}{n_{\Jdx}}
\newcommand{\sons}{\mathop{\operatorname{sons}}\nolimits}
\newcommand{\level}{\mathop{\operatorname{level}}\nolimits}
\newcommand{\treeroot}{\mathop{\operatorname{root}}\nolimits}
\newcommand{\pred}{\mathop{\operatorname{pred}}\nolimits}
\newcommand{\row}{\mathop{\operatorname{row}}\nolimits}

\newtheorem{theorem}{Theorem}
\newtheorem{definition}[theorem]{Definition}
\newtheorem{remark}[theorem]{Remark}
\begin{document}

\title{Computing the eigenvalues of symmetric ${\mathcal H}^2$-matrices
       by slicing the spectrum
\thanks{
    }
}

\author{Peter Benner
        \and Steffen B\"orm
        \and Thomas Mach
        \thanks{The research of the third author was partially supported
    by the Research Council KU Leuven, fellowship F+/13/020 Exploiting
    Unconventional QR-Algorithms for Fast and Accurate Computations of
    Roots of Polynomials; and by the Interuniversity Attraction Poles
    Programme, initiated by the Belgian State, Science Policy Office,
    Belgian Network DYSCO (Dynamical Systems, Control, and Optimization).}
       \and Knut Reimer
       \thanks{The research of the second and fourth authors was funded by
    the DFG Deutsche Forschungsgemeinschaft in project BO 3289/4-1.}}

\maketitle

\begin{abstract}
The computation of eigenvalues of large-scale matrices arising from
finite element discretizations has gained significant interest in the
last decade \cite{KMX13}.
Here we present an new algorithm based on slicing the spectrum that takes
advantage of the rank structure of resolvent matrices in order to compute
$m$ eigenvalues of the generalized symmetric eigenvalue problem in
${\mathcal O}(nm\log^\alpha n)$ operations, where $\alpha>0$ is a small
constant.

\end{abstract}

%
%
\section{Introduction}
\label{introduction}

The numerical solution of the generalized eigenproblem
\begin{equation}\label{eq:evp}
    \left( A - \lambda B \right) x = 0,
\end{equation}
given $A,B\in\mathbb{R}^{n\times n}$ and searching for $\lambda\in\mathbb{C}$
and $x\in\mathbb{C}^n\setminus\{0\}$, is one of the fundamental problems in
the computational sciences and engineering.
It arises in numerous applications ranging from structural and vibrational
analysis to problems in computational physics and chemistry like
electronic and band structure calculations, see, e.g., \cite{KMX13} and the
reports therein.
In particular, the investigation and design of new materials poses numerous
new challenges for the numerical solution of \eqref{eq:evp}.
These include the necessity to compute more than just the (few) smallest
magnitude eigenvalue(s) --- the target of many algorithms discussed in the
literature.
Often in these problems, a large number of interior eigenvalues are required.
This poses a significant challenge for most popular algorithms used to
solve large-scale eigenproblems based on the Arnoldi or Lanczos processes
or the Jacobi-Davidson method.
Therefore, we will discuss here a different approach that has received
little attention in the literature so far: the slicing-the-spectrum
approach discussed in \cite{Par80}.

Many of the application problems listed above lead to a symmetric eigenproblem
in the sense that $A=A^T$ and $B=B^T$.
Moreover, in applications arising from the discretization of (elliptic)
partial differential operators --- which probably cover the majority of
these application problems --- the matrix $B$ is a mass matrix and thus
positive definite, which we denote by $B>0$.
In this situation, the eigenvalues $\lambda$ and eigenvectors $x$ are all
real.
Here, we will assume these conditions and furthermore, we will focus on the
computation of eigenvalues.
If necessary, eigenvectors corresponding to selected eigenvalues can be
computed by inverse iteration which we will not further discuss.

Slicing-the-spectrum allows to compute a selected number of eigenvalues of
a symmetric matrix, or even all of them.
It requires knowledge of the inertia of shifted versions of the matrix, which
can be computed by the $LDL^T$ factorization.
As this is quite an expensive computation for sparse matrices, the method
has received little attention in the literature.
For data-sparse matrices which allow a low complexity computation of the
$LDL^T$ factorization, though, this method becomes attractive again.
In \cite{BenM11c}, we have used this approach to show that some, say $m$,
eigenvalues of $\mathcal{H}_{\ell}$-matrices can be computed in
$\mathcal{O}(mn\log^\alpha(n))$ complexity (for a discussion of the involved
constants we refer to \cite{BenM11c}).
$\mathcal{H}_{\ell}$-matrices are a class of simple hierarchical
($\mathcal{H}$-) matrices that are rank-structured in their off-diagonal
parts.
That is, the off-diagonal parts of these matrices are represented in a
hierarchical way by low-rank blocks so that the total storage for the
matrix is of linear-logarithmic complexity.
Such matrices often arise from the discretization of non-local operators
arising in integral equations or as solution operators of (elliptic) partial
differential operators \cite{BoeGH02,Hac99,BebH03,Boe07,Hac09}, and can
therefore often be used in the above application problems for the algebraic
representation of the involved integral and differential operators.

It was shown in \cite{BenM11c} that the $LDL^T$ factorization for
$\mathcal{H}_{\ell}$-matrices has bounded block ranks.
This allows the efficient implementation of the slicing-the-spectrum
approach for these special $\mathcal{H}$-matrices.
Numerical experiments however illustrated that this does not hold for
$\mathcal{H}$-matrices, casting doubt on the usefulness of this approach
for more general rank-structured matrices.
In this paper, we investigate the slicing-the-spectrum approach for
$\mathcal{H}^{2}$-matrices.
This matrix format allows a further compression compared to
$\mathcal{H}$-matrices by considering the low-rank structure of the whole
off-diagonal part of a block-row rather than of individual blocks.
We will see that an efficient $LDL^T$ factorization of
$\mathcal{H}^{2}$-matrices is possible and thus, an efficient implementation
of the slicing-the-spectrum approach is feasible.
We will furthermore extend this approach from the standard eigenvalue
problem considered in \cite{BenM11c} to the symmetric-definite eigenproblem
\eqref{eq:evp} with $A,B$ symmetric and $B>0$.
Moreover, this approach is shown to be easily parallelizable which allows to
gain further efficiency on current computer architectures.
 
The paper is structured as follows: in Section~\ref{sec:h2matrix}, we
introduce the necessary background on $\mathcal{H}^{2}$-matrices.
We then discuss the efficient implementation of the $LDL^T$ factorization
in the $\mathcal{H}^{2}$-format.
The slicing-the-spectrum approach is then reviewed in
Section~\ref{sec:slicing}.
Furthermore, the application to $\mathcal{H}^{2}$-matrices is discussed
as well as the extension to the symmetric-definite eigenproblem.
We also discuss a parallel implementation of the method.
Numerical experiments illustrating the performance of the
$\mathcal{H}^{2}$-slicing-the-spectrum algorithm and its parallelization are
presented in Section~\ref{sec:numerexp}.

%
%

\section{\texorpdfstring{$\mathcal{H}^{2}$}{H2}-Matrices and Their \texorpdfstring{$LDL^T$}{LDLT} Factorization}
\label{sec:h2matrix}

\subsection{\texorpdfstring{${\mathcal H}^2$}{H2}-matrices}

In this section we briefly recollect the basic definitions of
$\mathcal{H}^2$-matrices \cite{HAKHSA00,BO10}:
matrices are split into submatrices according to a \emph{block tree},
and this tree is constructed using \emph{cluster trees} describing
the decomposition of row and co\-lumn index sets.
If a submatrix is \emph{admissible}, it is represented in
factorized form using low-rank \emph{cluster bases} and
\emph{coupling matrices}.

%
%
\begin{definition}[Cluster tree]
\label{definition: cluster tree}
Let $\Idx$ be a finite index set, and let ${\mathcal T}$ be a
labeled tree.
Denote the label of each node $t\in{\mathcal T}$ by $\hat t$.

${\mathcal T}$ is called a \emph{cluster tree} for $\Idx$ if
the following conditions hold:
\begin{itemize}
  \item its root $r=\treeroot({\mathcal T})$ is labeled by $\Idx$,
        i.e., $\hat r=\Idx$,
  \item for all $t\in{\mathcal T}$ with $\sons(t)\neq\emptyset$,
        we have $\hat t = \bigcup_{t'\in\sons(t)} \hat t'$,
  \item for all $t\in{\mathcal T}$ and all $t_1,t_2\in\sons(t)$,
        $t_1\neq t_2$, we have $\hat t_1\cap\hat t_2=\emptyset$.
\end{itemize}
A cluster tree for $\Idx$ is denoted by $\ctI$, its nodes are called
\emph{clusters}, and $\lfI := \{ t\in\ctI\ :\ \sons(t)=\emptyset \}$
defines the set of its leaves.
\end{definition}

%
%
\begin{remark}[Leaf partition]
\label{remark: cluster tree 1}
The definition implies $\hat t\subseteq\Idx$ for all $t\in\ctI$.

We also have that the labels of the leaves of $\ctI$ form a disjoint
partition $\{ \hat t\ :\ t\in\lfI \}$ of the index set $\Idx$
\cite{HacK00,BO10}.
\end{remark}

%
%
\begin{remark}[Cardinalities]
\label{remark: cluster tree 2}
Let $\nI := \#\Idx$ denote the number of indices.
In typical situations, a cluster tree consists of ${\mathcal O}(\nI/k)$
clusters, where $k$ denotes the rank used to approximate matrix blocks.

The sum of the cardinalities of the index sets corresponding
to all clusters is typically in
${\mathcal O}(\nI \log(\nI))$ \cite{BO10}, since each index appears
in ${\mathcal O}(\log(\nI))$ clusters.
\end{remark}

Remarks \ref{remark: cluster tree 1} and \ref{remark: cluster tree 2} imply
that algorithms with optimal (linear) complexity should have at most constant
complexity in all non-leaf clusters $t\in\ctI\setminus\lfI$ and linear
complexity (with respect to $\#\hat t$) in all leaf clusters $t\in\lfI$.

With the help of the cluster tree we are able to define the block tree,
which gives us a hierarchically structured block partition of
$\Idx\times\Jdx$ and ultimately a partition of matrices
$X\in\bbbr^{\Idx\times\Jdx}$ into submatrices.

%
%
\begin{definition}[Block tree]
\label{definition: block tree}
Let $\ctI$ and $\ctJ$ be cluster trees for index sets $\Idx$
and $\Jdx$, respectively.

A labeled tree ${\mathcal T}$ is called a \emph{block tree} for
$\ctI$ and $\ctJ$ if the following conditions hold:
\begin{itemize}
  \item for all $b\in{\mathcal T}$, there are $t\in\ctI$
        and $s\in\ctJ$ such that $b=(t,s)$ and $\hat b=\hat t\times\hat s$,
  \item the root is
        $r=\treeroot({\mathcal T})=(\treeroot(\ctI),\treeroot(\ctJ))$,
  \item for all $b=(t,s)\in{\mathcal T}$, $\sons(b)\neq\emptyset$,
        we have
        \begin{equation*}
          \sons(b) = \begin{cases}
            \{t\}\times\sons(s)
            &\text{ if } \sons(t)=\emptyset\neq\sons(s),\\
            \sons(t)\times\{s\}
            &\text{ if } \sons(t)\neq\emptyset=\sons(s),\\
            \sons(t)\times\sons(s)
            &\text{ otherwise}.
          \end{cases}
        \end{equation*}
\end{itemize}
A block tree for $\ctI$ and $\ctJ$ is denoted by $\ctIJ$, its nodes
are called \emph{blocks}, and the set of its leaves is denoted by
$\lfIJ := \{ b\in\ctIJ\ :\ \sons(b)=\emptyset \}$.

For all blocks $b=(t,s)\in\ctIJ$, $t$ is called the \emph{row cluster}
 and $s$ is called the \emph{column cluster}.
\end{definition}

%
%
\begin{remark}[Leaf partition]
The Definitions~\ref{definition: cluster tree} and
\ref{definition: block tree} imply that a block tree $\ctIJ$ is a
cluster tree for $\Idx\times\Jdx$ and that therefore the set of
the labels of its leaves $\{ \hat t\times\hat s \ :\ b=(t,s)\in\lfIJ \}$
is a disjoint partition of $\Idx\times\Jdx$.
We use this partition to split matrices into submatrices.
\end{remark}

To determine which of these submatrices can be approximated
by low-rank representations, we split the set $\lfIJ$ of leaf blocks
into a set of \emph{admissible} blocks and a remainder of
``sufficiently small'' blocks.

%
%
\begin{definition}[Admissible blocks]
Let $\lfaIJ\subseteq\lfIJ$ be a subset of the leaves $\lfIJ$ and let
$\lfiIJ := \lfIJ\setminus\lfaIJ$ denote the remaining leaves.

If $(t,s)\in\lfiIJ$ implies $t\in\lfI$ and $s\in\lfJ$, we call
$\lfaIJ$ a set of \emph{admissible blocks} and $\lfiIJ$ the
corresponding set of \emph{inadmissible blocks}.
\end{definition}

Typically we choose the set $\lfaIJ$ of admissible leaves in
a way that ensures that for each $b=(t,s)\in\lfaIJ$, the corresponding
submatrix $X_{|\hat t\times\hat s}$ can be approximated by low rank.
In practice a minimal block tree is constructed based on an
\emph{admissibility condition} that predicts whether a given block
$b=(t,s)$ can be approximated.
If this is the case, the block is chosen as an admissible leaf of $\ctIJ$.
Otherwise we either check the sons of $b$ given by
Definition~\ref{definition: block tree} or, if there are no sons,
declare the block an inadmissible leaf.

In the context of elliptic partial differential equations, we usually
employ an \emph{admissibility criterion} of the form
\begin{equation*}
  \max\{\operatorname{diam}(t), \operatorname{diam}(s)\}
  \leq 2\eta \operatorname{dist}(t, s),
\end{equation*}
where $\operatorname{diam}(t)$ and $\operatorname{dist}(t,s)$ denote the
diameter and distance of clusters in a suitable way.

%
%
\begin{remark}[Sparse block tree]
\label{remark: block tree}
If there is a constant $c_{\rm sp}\in\bbbn$ such that
\begin{align*}
  \#\{ s\in\ctJ\ :\ (t,s)\in\ctIJ \} &\leq c_{\rm sp} &
  &\text{ for all } t\in\ctI,\\
  \#\{ t\in\ctI\ :\ (t,s)\in\ctIJ \} &\leq c_{\rm sp} &
  &\text{ for all } s\in\ctJ
\end{align*}
hold, we call the block tree $\ctIJ$ \emph{$c_{\rm sp}$-sparse}.

In this case, Remark~\ref{remark: cluster tree 2} implies that the number
of blocks $\#\ctIJ$ is in ${\mathcal O}(\nI+\nJ)$, so algorithms of
optimal complexity should require only a constant number of operations
per block.
\end{remark}

$\mathcal H^2$-matrices use a three-term representation $V_t S_{(t,s)} W_s^T$
for all admissible blocks.
The matrix $V_t$ depends only on the row cluster $t$ and $W_s$ only on the
column cluster $s$.
The advantage of this representation is that only the $k\times k$ matrix
$S_{(t,s)}$is stored for every admissible block $(t,s)$.

Storing the matrices $V_t$ and $W_s$ directly would lead to linear complexity
in each cluster.
Thus we would get log-linear complexity for the whole families of matrices
$(V_t)_{t\in\ctI}$ and $(W_s)_{s\in\ctJ}$ (see Remark~\ref{remark: cluster tree 2}).
In \cite{HacKS00} the more efficient \emph{nested} representation of these
families is introduced.

%
%
\begin{definition}[Cluster basis]
Let $k\in\bbbn$, and let $(V_t)_{t\in\ctI}$ be a family of matrices
satisfying $V_t\in\bbbr^{\hat t\times k}$ for all $t\in\ctI$.

This family is called a \emph{(nested) cluster basis} if for
each $t\in\ctI$ there is a matrix $E_t\in\bbbr^{k\times k}$ such
that
\begin{align}\label{eq:nested}
  V_{t|\hat t'\times k} &= V_{t'} E_{t'} &
  &\text{ for all } t\in\ctI,\ t'\in\sons(t).
\end{align}
The matrices $E_t$ are called \emph{transfer matrices}, and
$k$ is called the \emph{rank} of the cluster basis.
\end{definition}

Due to (\ref{eq:nested}), we only have to store the $\hat t\times k$ matrices
$V_t$ for leaf clusters $t\in\lfI$ and the $k\times k$ transfer matrices
$E_t$ for all clusters $t\in\ctI$.

%
%
\begin{remark}[Storage]
\label{remark:storage cluster basis}
According to Remark~\ref{remark: cluster tree 1}, the ``leaf matrices''
$(V_t)_{t\in\lfI}$ require $\nI k$ units of storage.
The transfer matrices $(E_t)_{t\in\ctI}$ require $k^2 \#\ctI$ units of
storage.
With the standard assumption $\#\ctI\lesssim \nI/k$, we can conclude
that a cluster basis can be represented in ${\mathcal O}(\nI k)$ units
of storage \cite{HacKS00,BOHA02,BO10}.
\end{remark}

%
%
\begin{definition}[${\mathcal H}^2$-matrix]
Let $\ctI$ and $\ctJ$ be cluster trees for index sets $\Idx$ and
$\Jdx$, let $\ctIJ$ be a matching block tree, and let
$(V_t)_{t\in\ctI}$ and $(W_s)_{s\in\ctJ}$ be nested cluster bases.

A matrix $G\in\bbbr^{\Idx\times\Jdx}$ is called an
\emph{${\mathcal H}^2$-matrix} for $\ctIJ$, $(V_t)_{t\in\ctI}$ and
$(W_s)_{s\in\ctJ}$ if for each admissible block $b=(t,s)\in\lfaIJ$
there is a matrix $S_b\in\bbbr^{k\times k}$ such that
\begin{equation}\label{eq:vsw}
  G_{|\hat t\times\hat s} = V_t S_b W_s^T.
\end{equation}
The matrices $S_b$ are called \emph{coupling matrices}, the cluster
bases $(V_t)_{t\in\ctI}$ and $(W_s)_{s\in\ctJ}$ are called \emph{row} and
\emph{column cluster bases}.
\end{definition}

%
%
\begin{remark}[Storage]
An $\mathcal H^2$-matrix is represented by its nested cluster bases,
its $k\times k$ coupling matrices $(S_b)_{b\in\lfIJ^+}$ and its nearfield
matrices $(G_{|\hat t\times\hat s})_{b\in\lfIJ^-}$.
We have already seen in Remark~\ref{remark:storage cluster basis} that
the nested representations of the cluster bases require
${\mathcal O}(\nI k)$ and ${\mathcal O}(\nJ k)$ units of storage,
respectively.
The coupling matrices require ${\mathcal O}(k^2)$ units of storage
per block, leading to total requirements of ${\mathcal O}(\nI k)$
for a sparse block tree $\ctIJ$.
For $(t,s)=b\in\lfIJ^-$ both $t$ and $s$ are leaf clusters and so
$\#\hat t$ and $\#\hat s$ are small, usually bounded by $k$, and we
can conclude that the nearfield matrices require ${\mathcal O}(\nI k)$
units of storage if $\ctIJ$ is sparse.
In total an ${\mathcal H}^2$-matrix representation requires only
${\mathcal O}((\nI+\nJ) k)$ units of storage \cite{BOHA02,BO10}.
\end{remark}

Approximating an arbitrary matrix $X\in\bbbr^{\Idx\times\Jdx}$ by an
${\mathcal H}^2$-matrix becomes a relatively simple task if we
apply orthogonal projections.
These projections are readily available if the cluster bases
are \emph{orthogonal}:

%
%
\begin{definition}[Orthogonal cluster basis]
We call a cluster basis $(V_t)_{t\in\ctI}$ \emph{orthogonal} if
\begin{align*}
  V_t^T V_t &= I &
  &\text{ for all } t\in\ctI.
\end{align*}
\end{definition}

If $(V_t)_{t\in\ctI}$ and $(W_s)_{s\in\ctJ}$ are orthogonal cluster bases, the
optimal coupling matrices for a given matrix $G$ (with respect both to the
Frobenius norm and the spectral norm) are given by
\begin{align}
\label{ali: optimal coupling}
  S_b &:= V_t^T G_{|\hat t\times\hat s} W_s &
  &\text{ for all } b=(t,s)\in\lfaIJ.
\end{align}
This property can be used to compute the best approximation of the
product of ${\mathcal H}^2$-matrices in ${\mathcal O}(n k^2)$ operations
\cite{BO04a} as long as both cluster bases are known in advance.
Unfortunately the suitable cluster bases for the results of arithmetic
operations are typically not known.
Thus we have to construct adaptive cluster bases during the computations, see section \ref{subsec: Adaptive Cluster Basis} and \cite{BO05a,BO10,BOHA02}.

\subsection{Algebraic operations}

We want to compute the eigenvalues of a matrix $A\in\bbbr^{\Idx\times\Idx}$
corresponding to a Galerkin discretization of an elliptic partial
differential equation via a slicing method.
This method relies on a sufficiently accurate approximation of the
$LDL^T$ factorization of shifted matrices.

In order to construct an approximation of this factorization, we employ
an algorithm based on low-rank updates \cite{BORE14}.
We assume for the sake of simplicity that every non-leaf cluster has exactly two sons.
We obtain the following block equation for the $LDL^T$ factorization of
a submatrix $A_{|\hat t\times\hat t}$ for non-leaf clusters $t$ with
$\sons(t)=\{t_1,t_2\}$:
\begin{align*}
  \begin{pmatrix}
    A_{11} & A_{12} \\ A_{21} & A_{22}
  \end{pmatrix}
  \ &=\ A_{\hat{t}\times\hat{t}} = L_{\hat{t}\times\hat{t}}D_{\hat{t}\times\hat{t}}L_{\hat{t}\times\hat{t}}^T\\
  \ &=\ \begin{pmatrix}
         L_{11} &  \\ L_{21} & L_{22}
       \end{pmatrix}
       \begin{pmatrix}
         D_{11} &  \\ & D_{22}
       \end{pmatrix}
       \begin{pmatrix}
         L_{11} & L_{21}^T \\  & L_{22}
       \end{pmatrix}
  \\ &=\ \begin{pmatrix}
         L_{11}D_{11}L_{11}^T & L_{11}D_{11}L_{21}^T \\ L_{21}D_{11}L_{11}^T & L_{21}D_{11}L_{21}^T + L_{22}D_{22}L_{22}^T
       \end{pmatrix}.
\end{align*}
We can solve $A_{11} = L_{11}D_{11}L_{11}^T$ by recursion to get $L_{11}$
and $D_{11}$.
If the recursion reaches a leaf block, the block is a sufficiently small
matrix in standard representation and the $LDL^T$ factorization can be
computed by standard algorithms.

In the second step we can obtain $L_{21}$ by solving the triangular system
$A_{12} = L_{11}D_{11}L_{21}^T$.
This requires forward substitution for $A_{12} = L_{11}Y$ and solving the
diagonal system $Y = D_{11}L_{21}^T$.
The same block equation approach as above reduces the forward substitution
to matrix-matrix multiplications of the form $C \gets C + \alpha A B$.

Finally we can solve $A_{22} - L_{21}D_{11}L_{21}^T = L_{22}D_{22}L_{22}^T$ to
get $L_{22}$ and $D_{22}$.
This means a matrix-matrix multiplication of the form $C \gets C + \alpha A B$
and a recursion as in step one.

The block equation approach for the matrix-matrix multiplication
$C \gets C + \alpha A B$ leads to recursive calls
$C_{ij} \gets C_{ij} + \alpha A_{ik} B_{kj}$.
The basis case of the recursion is when $A$ or $B$ is a leaf.
Admissible leaves have low rank because of their three-term representation.
Inadmissible leaves have low rank because they are small.
In both cases we can compute a low rank representation $X Y^T$ of the product
$A B$ in linear complexity.

Altogether the arithmetic is reduced to the task of applying low-rank updates
$C_{|\hat t\times\hat r} + X Y^T$ to a submatrix of an ${\mathcal H}^2$-matrix,
where $X\in\bbbr^{\hat t\times k}$ and $Y\in\bbbr^{\hat s\times k}$.

\subsection{${\mathcal H}^2$-matrix Representation of $C + X Y^T$}

In order to handle low-rank updates to $\mathcal{H}^2$-matrices
efficiently, we follow the approach described in \cite{BORE14}, i.e.,
we consider $C + X Y^T$ as an $\mathcal{H}^2$-matrix with increased
rank and apply the recompression algorithm \cite{BOHA02} in order
to reduce the rank while guaranteeing a given accuracy.
We only outline the algorithm here for the sake of completeness
and refer readers to \cite{BORE14} for details.

We first consider a \emph{global} low-rank update $C \gets C + X Y^T$
and start by examining the ${\mathcal H}^2$-matrix representation of
the new matrix $\widetilde{C}:=C + X Y^T$.
For each admissible leaf $b=(t,s)\in\lfaIJ$, we obtain the following
simple three-term representation:
\begin{align*}
  \widetilde{C}_{|\hat t\times\hat s}
  &= (C+X Y^T)_{|\hat t\times\hat s}
   = V_t S_b W_s^T + X_{|\hat t\times k} Y_{|\hat s\times k}^T \\
  &= \begin{pmatrix} V_t & X_{|\hat t\times k} \end{pmatrix}
    \begin{pmatrix} S_b & \\ & I_k \end{pmatrix}
    \begin{pmatrix} W_s & Y_{|\hat s\times k} \end{pmatrix}^T.
\end{align*}
This leads to the new cluster bases
\begin{align*}
  \widetilde V_t = \begin{pmatrix} V_t & X_{|\hat t\times k} \end{pmatrix}
  \quad \text{ and } \quad 
  \widetilde W_s = \begin{pmatrix} W_s & Y_{|\hat s\times k} \end{pmatrix}.
\end{align*}
These are nested with transfer matrices
\begin{align*}
  \widetilde E_t = \begin{pmatrix} E_t & \\ & I_k \end{pmatrix}
  \quad \text{ and } \quad 
  \widetilde F_s = \begin{pmatrix} F_s & \\ & I_k \end{pmatrix}.
\end{align*}
The new nested cluster bases $\widetilde V$ and $\widetilde W$ together
with coupling matrices
\begin{equation*}
  \widetilde S_b = \begin{pmatrix} S_{b} & \\ & I_k \end{pmatrix}
\end{equation*}
for each $b\in\lfIJ$ give us an \emph{exact} ${\mathcal H}^2$-matrix
representation of $\widetilde{C}=C + XY^T$.

The drawback of this representation is the doubled rank.
We solve this problem by applying the recompression algorithm described
in \cite{BO05a,BO10}:
we construct adaptive orthogonal cluster bases and then approximate the
original matrix in the space defined by these bases
(cf. (\ref{ali: optimal coupling})).

\subsection{Weight Matrices}

In order to keep the presentation simple we denote the $\mathcal{H}^2$-matrix
$\widetilde{C}$ in the following just by $C$, the corresponding row and
column cluster bases by $(V_t)_{t\in\ctI}$ and $(W_s)_{s\in\ctJ}$, their rank by $k$, the
coupling matrices by $(S_b)_{b\in\lfaIJ}$, and the nearfield matrices
by $(C_{|\hat b})_{b\in\lfiIJ}$.
In our algorithm, these matrices are constructed implicitly according to
the equations given in the previous section.

The recompression algorithm is based on the method introduced
in \cite{BOHA02} using the refinements added in \cite{BO05a}:
the original algorithm relies on approximations of certain
submatrices of $C$, and since this is an $\mathcal{H}^2$-matrix,
these submatrices can be represented by compact \emph{weight matrices}.
Here we only briefly outline the concept and refer readers
to \cite{BO05a} and \cite[Chapter~6.6]{BO10} for details.

We consider only the construction of a row basis, since a column
basis can be obtained by applying the same algorithm to the
transposed matrix $C^T$.

The cluster basis $V_t$ is directly used for the representation of all
admissible blocks $(t,s)\in\lfaIJ$.
We collect the corresponding column clusters in the set
\begin{equation*}
  \row(t) = \{s\in\ctJ | (t,s)\in\lfaIJ\}.
\end{equation*}
Because of the nested structure, $V_t$ influences also blocks
$(t^*,s)\in\lfaIJ$ connected to predecessors $t^*$ of $t$.
We denote the set of predecessors by
\begin{equation*}
  \pred(t) := \begin{cases}
    \{t\} &\text{ if } t=\treeroot(\ctI),\\
    \{t\} \cup \pred(t^+)
    &\text{ for } t^+\in\ctI,t\in\sons(t^+).
  \end{cases}
\end{equation*}
For the construction of the new cluster basis $(Q_t)_{t\in\ctI}$, we have to
consider the set
\begin{equation*}
  \row^*(t) = \bigcup_{t^*\in\pred(t)}\row(t^*).
\end{equation*}
Let $\row(t) = \{s_1,...,s_\sigma\}$ and $\row^*(t) = \{s_1,...,s_\rho\}$
with $\sigma=\#\row(t)$ and $\rho=\#\row^*(t)$.
The part of $C$ which is described (directly or indirectly) by $V_t$ is
\begin{equation*}
  C_t := \begin{pmatrix} C_{|\hat t\times\hat s_1} & \ldots
                & C_{|\hat t\times\hat s_\rho} \end{pmatrix}.
\end{equation*}
Using the approach of (\ref{ali: optimal coupling}) we search for an
orthogonal matrix $Q_t$ with lower rank than $V_t$ such that 
\begin{equation*}
  Q_tQ_t^T C_t \approx C_t.
\end{equation*}
We want to reach this goal via singular value decomposition.

Computing the SVD of $C_t$ directly would be too expensive, but we
can introduce weight matrices to significantly reduce the number of
operations:
if for a matrix $Z_t$ there is an orthogonal matrix $P_t$ with
\begin{equation*}
  C_t = V_t Z_t^T P_t^T,
\end{equation*}
we call $Z_t$ a weight matrix for $C$ and $t$.
Since $V_t$ has only $k$ columns, the same holds for $Z_t$, and using,
e.g., a QR factorization leads to an upper triangular $Z_t$ with $k$ columns
and not more than $k$ rows.
For the construction of the cluster basis, we are only interested in
the left singular vectors and the singular values of $C_t$.
Due to the orthogonality of $P_t$, these quantities can be obtained by
computing only the SVD of $V_t Z_t^T$ instead of working with $C_t$.
Since $Z_t^T$ has not more than $k$ columns, this approach leads to
a significant reduction in the computational work.

We construct the weight matrices by a top-down recursion:
for the root of $\ctI$, the weight matrix can be computed directly.
For a cluster $t\in\ctI\setminus\{\treeroot(\ctI)\}$, we assume
that a weight matrix $Z_{t^+}$ for its father $t^+\in\ctI$ has already
been computed and denote the corresponding orthogonal matrix by
$P_{t^+}$.
Since the cluster basis $(V_t)_{t\in\ctI}$ is nested, we have
\begin{align}
  C_t
  &= \begin{pmatrix}
       C_{|\hat t\times\hat s_1} & \ldots & C_{|\hat t\times\hat s_\sigma}
       & C_{t^+|\hat t\times\Jdx}
     \end{pmatrix}\notag\\
  &= \begin{pmatrix}
       C_{|\hat t\times\hat s_1} & \ldots & C_{|\hat t\times\hat s_\sigma}
       & (V_{t^+} Z_{t^+}^T P_{t^+}^T)_{|\hat t\times\Jdx}
     \end{pmatrix}\notag\\
  &= \begin{pmatrix}
       C_{|\hat t\times\hat s_1} & \ldots & C_{|\hat t\times\hat s_\sigma}
       & V_{t^+|\hat t\times k} Z_{t^+}^T P_{t^+}^T
     \end{pmatrix}\notag\\
  &= \begin{pmatrix}
       V_t S_{(t,s_1)} W_{s_1}^T & \ldots
       & V_t S_{(t,s_\sigma)} W_{s_\sigma}^T & V_t E_t Z_{t^+}^T P_{t^+}^T
     \end{pmatrix}
   = V_t B_t\label{equ:VtBt}
\end{align}
with the matrix
\begin{equation*}
  B_t := \begin{pmatrix}
    S_{(t,s_1)} W_{s_1}^T &
    \ldots &
    S_{(t,s_\sigma)} W_{s_\sigma}^T &
     E_t Z_{t^+}^T P_{t^+}^T
  \end{pmatrix}.
\end{equation*}
This allows us to obtain the following factorized representation
of $C_t$:
\begin{align}
  C_t 
  &= V_t B_t
   = V_t \begin{pmatrix}
       S_{(t,s_1)} W_{s_1}^T & \cdots
       & S_{(t,s_\sigma)} W_{s_\sigma}^T
       & E_t Z_{t^+}^T P_{t^+}^T
     \end{pmatrix}\notag\\
  &= V_t
     \begin{pmatrix}
       S_{(t,s_1)} & \cdots & S_{(t,s_\sigma)} & E_t Z_{t^+}^T
     \end{pmatrix}
     \begin{pmatrix}
       W_{s_1} & & & \\
       & \ddots & & \\
       & & W_{s_\sigma} & \\
       & & & P_{t^+}
     \end{pmatrix}^T\notag\\
  &= V_t\widetilde Z_t^T \widetilde P_t^T.\label{equ: weight matrix}
\end{align}
We assume in the following that the cluster basis $W$ is orthogonal.
If it is not, we can apply recursive QR factorizations to replace it
by an orthogonal basis in linear complexity \cite[Section~3.2]{Boe09}.
Then $\widetilde P_t$ is orthogonal and $\widetilde Z_t$ is a weight
matrix, but the number of rows of $\widetilde Z_t$ typically exceeds
the number of columns.
Thus we compute a thin QR decomposition $\widetilde Z_t = \hat P_t Z_t$ and
get
\begin{equation*}
  C_t = V_t \widetilde Z_t^T \widetilde P_t^T
  = V_t Z_t^T \hat P_t^T \widetilde P_t^T
  = V_t Z_t^T P_t^T.
\end{equation*}
$P_t$ is orthogonal, and so $Z_t$ is a small $k\times k$ weight matrix.

Altogether we can compute the weight matrices by a top down algorithm
which only assembles $\widetilde Z_t$ and computes its QR decomposition.
Only $k\times k$ weight matrices $Z_t$ are stored and the number of
considered blocks $\sigma$ is bounded by the constant $c_{sp}$.
Thus the storage requirement for one cluster $t\in\ctI$ is in
${\mathcal O}(k^2)$ and the computational time is in ${\mathcal O}(k^3)$.
The storage requirement for all weight matrices is in
${\mathcal O}(k^2\#\ctI)$ and the computational time for the whole
algorithm is in ${\mathcal O}(k^3\#\ctI)$ \cite{BO05a,BO10}.
Using the standard assumption $\#\ctI\lesssim \nI/k$, we conclude that
${\mathcal O}(\nI k)$ units of storage and ${\mathcal O}(\nI k^2)$
operations are sufficient to set up all weight matrices.

\subsection{Adaptive Cluster Basis}
\label{subsec: Adaptive Cluster Basis}

The weight matrices can be computed efficiently by a top-down traversal
of the cluster tree $\ctI$.
Once they are at our disposal, we can use a bottom-up traversal of the
cluster tree to construct the required adaptive cluster basis
$(Q_t)_{t\in\ctI}$ following the method given in \cite{BO05a} and
\cite[Chapter~6.6]{BO10}.

With the help of the weight matrices we get
\begin{equation}\begin{split}
\label{ali: error 1}
  \|Q_t Q_t^T C_t - C_t\| 
  &= \|Q_t Q_t^T V_t Z_t^T P_t^T - V_t Z_t^T P_t^T\| \\
  &= \ \|Q_t Q_t^T V_t Z_t^T - V_t Z_t^T\|
\end{split}\end{equation}
for both the spectral and the Frobenius norm.
Thus we only have to compute the SVD of $V_tZ_t^T$ instead of $C_t$.
The direct approach would have linear complexity in each cluster and 
we would end up with log-linear complexity due to
Remark~\ref{remark: cluster tree 2}.
We also would not obtain a nested cluster basis.

In order to avoid both problems, we take advantage of the nested
structure of $(V_t)_{t\in\ctI}$ and $(Q_t)_{t\in\ctI}$.
We arrange the computation of the adaptive cluster basis $(Q_t)_{t\in\ctI}$
in a bottom-up algorithm that also computes the basis change matrices
$R_t := Q_t^T V_t$ for all $t\in\ctI$ that can be used to compute the
new coupling matrices efficiently.

In leaf clusters we compute the SVD of $V_tZ_t^T$ directly and use the
left singular vectors corresponding to the $k$ largest singular values
to construct the orthogonal matrix $Q_t$.
The computational time for each leaf is $O(k^2\#\hat t)$ and for all
leaves together it is in $O(\nI k^2)$
(see Remark~\ref{remark: cluster tree 1}). 

The cluster basis in a non-leaf cluster is given by the nested
representation
\begin{equation*}
  V_t = \begin{pmatrix} V_{t_1} E_{t_1}\\
                        V_{t_2} E_{t_2} \end{pmatrix}.
\end{equation*}
We assume that the matrices $Q_{t_1}$ and $Q_{t_2}$ for the sons
have already been computed, and the nested structure of $(Q_t)_{t\in\ctI}$
implies that anything that cannot be represented by these matrices
also cannot be represented by $Q_t$, so applying a projection to
the range of the son matrices does not change the quality of the
approximation.
If we let
\begin{equation*}
  U_t := \begin{pmatrix}
    Q_{t_1} & \\
    & Q_{t_2}
  \end{pmatrix},
\end{equation*}
the orthogonal projection is given by $U_t U_t^T$ and applying it to
$V_t$ yields
\begin{align*}
  U_t U_t^T V_t
  = U_t
    \begin{pmatrix} Q_{t_1}^T & \\ & Q_{t_2}^T \end{pmatrix}
    \begin{pmatrix} V_{t_1}E_{t_1} \\ V_{t_2}E_{t_2}\end{pmatrix}
  = U_t
    \begin{pmatrix} R_{t_1}E_{t_1} \\ R_{t_2}E_{t_2}\end{pmatrix} 
  = U_t \widehat V_t
\end{align*}
with a $(2k)\times k$ matrix $\widehat V_t = U_t^T V_t$.
We compute the SVD of $\widehat V_t Z_t$ and again use the left singular vectors corresponding
to the $k$ largest singular values to form an orthogonal matrix
$\widehat Q_t\in\bbbr^{(2k)\times k}$.
The new cluster basis is defined by $Q_t := U_t \widehat Q_t$.
We deduce with Pythagoras' identity
\begin{align}
  \|Q_t &Q_t^T V_t Z_t^T - V_t Z_t^T\|^2 \notag\\
  &= \|U_t\widehat Q_t\widehat Q_t^TU_t^TV_tZ_t^T - U_tU_t^TV_tZ_t^T\|^2\notag\\
  &\quad + \|U_tU_t^TV_tZ_t^T - V_tZ_t^T\|^2\notag\\
  &= \|\widehat Q_t\widehat Q_t^T\widehat V_tZ_t^T - \widehat V_tZ_t^T\|^2 
  + \|U_tU_t^TV_tZ_t^T - V_tZ_t^T\|^2\label{equ: error 2}.
\end{align}
Thus the error for the cluster $t$ can be bounded by the error of the
projection of the son clusters and the error of the truncated SVD of
$\hat{V}_tZ_t^T$.
We will investigate the error in subsection \ref{subsec: error control}.

The basis change matrix $R_t$ is computed in ${\mathcal O}(k^3)$ operations
via
\begin{equation*}
  R_t = Q_t^T V_t = \widehat Q_t^T U_t^T V_t = \widehat Q_t^T \widehat V_t.
\end{equation*}
The transfer matrices of $Q_t$ can be constructed by splitting
$\widehat Q_t$ into its lower and upper half, i.e., by using
\begin{equation*}
  Q_t = U_t \widehat Q_t
  = \begin{pmatrix} Q_{t_1} & \\ & Q_{t_2} \end{pmatrix}
    \begin{pmatrix} F_{t_1} \\ F_{t_2} \end{pmatrix}.
\end{equation*}
We can see that leaf clusters $t\in\lfI$ require
${\mathcal O}(k^2 \#\hat t)$ operations while non-leaf clusters
$t\in\ctI\setminus\lfI$ require ${\mathcal O}(k^3)$.
The total computational time of the algorithm therefore is in
${\mathcal O}(k^2 \nI + k^3 \#\ctI)$ due to
Remark~\ref{remark: cluster tree 1}.
By the standard assumption $\#\ctI\lesssim\nI/k$, we conclude that not
more than ${\mathcal O}(\nI k^2)$ operations are required to construct
the new cluster basis \cite{BO05a,BO10}.

\subsection{Error Control}
\label{subsec: error control}

As we have seen in the previous subsections we are able to recompress
an $\mathcal{H}^2$-matrix in linear complexity and (\ref{equ: error 2})
suggests that the resulting error can be controlled by the accuracy of
the truncated SVD.
In this section, we describe a simplified version of the blockwise
error control strategy developed in \cite{BO05a} that, according to
our experiments, is suitable for treating eigenvalue problems.

Let $b=(t,s)\in\lfaIJ$.
Multiplying the matrices in (\ref{equ: error 2}) by $P_t^T$ from the
right, using $C_t = V_t B_t = V_t Z_t^T P_t^T$, and restricting to
$\hat t\times\hat s$, we obtain
\begin{align}
  \|Q_t &Q_t^T C_{|\hat t\times\hat s} - C_{|\hat t\times\hat s}\|^2\notag\\
  &= \|Q_t Q_t^T V_t B_{t|k\times\hat s}
         - V_t B_{t|k\times\hat s}\|^2\notag\\
  &= \|\widehat Q_t \widehat Q_t^T \widehat V_t B_{t|k\times\hat s}
      - \widehat V_t B_{t|k\times\hat s}\|^2\notag\\
  &\quad + \|U_t U_t^T V_t B_{t|k\times\hat s}
      - V_t B_{t|k\times\hat s}\|^2.\label{equ: error 2a}
\end{align}
Due to the nested structure of $V_t$ and the definition of $U_t$,
we have
\begin{align*}
  \|U_t &U_t^T V_t B_{t|k\times\hat s}
       - V_t B_{t|k\times\hat s}\|^2\\
  &= \|Q_{t_1} Q_{t_1}^T V_{t_1} E_{t_1} B_{t|k\times\hat s}
        - V_{t_1} E_{t_1} B_{t|k\times\hat s}\|^2\\
  &\quad +\|Q_{t_2} Q_{t_2}^T V_{t_2} E_{t_2} B_{t|k\times\hat s}
        - V_{t_2} E_{t_2} B_{t|k\times\hat s}\|^2\\
  &= \|Q_{t_1} Q_{t_1}^T V_{t_1} B_{t_1|k\times\hat s}
        - V_{t_1} B_{t_1|k\times\hat s}\|^2 \\
  &\quad +\|Q_{t_2} Q_{t_2}^T V_{t_2} B_{t_2|k\times\hat s}
        - V_{t_2} B_{t_2|k\times\hat s}\|^2.
\end{align*}
By simple induction we get
\begin{align}
  \|Q_t &Q_t^T C_{|\hat t\times\hat s} - C_{|\hat t\times\hat s}\|^2\notag\\ 
  &= \sum_{r\in\sons^*(t)}
    \|\widehat Q_r \widehat Q_r^T \widehat V_r B_{r|k\times\hat s}
        - \widehat V_r B_{r|k\times\hat s}\|^2\label{equ: error 3}
\end{align}
with the set of descendants given by
\begin{equation*}
  \sons^*(t)
  := \begin{cases}
       \{t\} & \text{if } t\in\lfI \\
       \{t\} \cup \bigcup_{t'\in\sons(t)} \sons^*(t') & \text{otherwise}
     \end{cases}
\end{equation*}
and extending the notation to $\widehat Q_t=Q_t$ and 
$\widehat V_t=V_t$ for leaf clusters $t\in\lfI$.

Equation (\ref{equ: error 3}) provides us with an explicit error
representation.
We get an efficiently computable error bound by extending
$B_{r|k\times\hat s}$ to the larger matrix $B_r$ and using
$B_r = Z_t^T P_t^T$ to reduce to the weight matrix:
\begin{align*}
  \|Q_t &Q_t^T C_{|\hat t\times\hat s} - C_{|\hat t\times\hat s}\|^2\\
  &= \sum_{r\in\sons^*(t)}
      \|\widehat Q_r \widehat Q_r^T \widehat V_r B_{r|k\times\hat s}
           - \widehat V_r B_{r|k\times\hat s}\|^2\\
  &\leq \sum_{r\in\sons^*(t)}
      \|\widehat Q_r \widehat Q_r^T \widehat V_r B_r
           - \widehat V_r B_r\|^2\\
  &= \sum_{r\in\sons^*(t)}
      \|\widehat Q_r \widehat Q_r^T \widehat V_r Z_r^T
           - \widehat V_r Z_r^T\|^2.
\end{align*}
This is an error bound that we can control directly via the truncation
criterion of the SVD used to compute $\widehat Q_r$.
Unfortunately it does not give us direct error  control for individual
blocks, which is crucial for efficient and reliable algebraic operations.
If we could bound each term in (\ref{equ: error 3}) by 
\begin{align*}
  \|\widehat Q_r &\widehat Q_r^T \widehat V_r B_{r|k\times\hat s}
     - \widehat V_r B_{r|k\times\hat s}\|^2\\
  &\leq \frac{\epsilon^2}{3} \|C_{|\hat t\times\hat s}\|^2
        \left(\frac{1}{3}\right)^{\operatorname{level}(r)-\operatorname{level}(t)},
\end{align*}
we would get
\begin{align}
  \|Q_t &Q_t^T C_{|\hat t\times\hat s} - C_{|\hat t\times\hat s}\|^2\notag\\
  &\leq \frac{\epsilon^2}{3} \|C_{|\hat t\times\hat s}\|^2
        \sum_{r\in\sons^*(t)}
           \left(\frac{1}{3}\right)^{\level(r)-\level(t)}\notag\\
  &= \frac{\epsilon^2}{3} \|C_{|\hat t\times\hat s}\|^2
        \sum_{\ell=\level(t)}^{p_\Idx}
           \left(\frac{1}{3}\right)^{\ell-\level(t)}
           \#\{ r\in\sons^*(t)\ :\ \level(r)=\ell \}\notag\\
  &\leq \frac{\epsilon^2}{3} \|C_{|\hat t\times\hat s}\|^2
        \sum_{\ell=\level(t)}^{p_\Idx}
           \left(\frac{2}{3}\right)^{\ell-\level(t)}
   < \epsilon^2 \|C_{|\hat t\times\hat s}\|^2\label{equ: error 4}
\end{align}
by the geometric summation formula.

We cannot simply set the tolerance in each cluster $r\in\sons^*(t)$ to
\begin{equation*}
  \omega_{r,b}^2
  := \frac{\epsilon^2}{3} \|C_{|\hat t\times\hat s}\|^2
        \left(\frac{1}{3}\right)^{\level(r)-\level(t)}
\end{equation*}
because it depends not only on $r$, but also on $b$.
The solution is to put the factor $\omega_{r,b}$ into the weight matrix
\cite{BO05a}.
The condition 
\begin{equation*}
  \|\widehat Q_r \widehat Q_r^T \widehat V_r B_{r|k\times\hat s}
       - \widehat V_r B_{r|k\times\hat s}\|^2
  \leq \omega_{r,b}^2
\end{equation*}
is equivalent to
\begin{equation*}
  \|\widehat Q_r \widehat Q_r^T \widehat V_r \omega_{r,b}^{-1}
          B_{r|k\times\hat s}
    - \widehat V_r \omega_{r,b}^{-1} B_{r|k\times\hat s}\|^2 \leq 1.
\end{equation*}
Since $\omega_{r',b} = \omega_{r,b}/3$ holds for all $r\in\sons^*(t)$
and $r'\in\sons(r)$, we can include the factors in the algorithm
for constructing the weight matrices in (\ref{equ: weight matrix}) and
get
\begin{equation}\label{equ: weighted weight matrix}
  \widetilde Z_{t,\omega}^T =
    \begin{pmatrix}
       \omega_{r,(t,s_1)}^{-1} S_{(t,s_1)}
       & \cdots
       & \omega_{r,(t,s_\sigma)}^{-1} S_{(t,s_\sigma)}
       & 3 E_t Z_{t^+}^T
    \end{pmatrix}.
\end{equation}
The resulting weight matrices $Z_{t,\omega}$ satisfy
\begin{equation*}
  (Z_{t,\omega}^T P_t^T)_{|k\times\hat s}
  = \omega_{r,b}^{-1} B_{r|k\times\hat s},
\end{equation*}
therefore we get the error bound in (\ref{equ: error 4}) if we replace
$Z_t$ by $Z_{t,\omega}$ and ensure that the rank $k$ used in the truncation
is large enough to capture all singular values larger than one.

Now we have found a recompression algorithm with linear complexity
${\mathcal O}(\nI k^2)$ allowing us to control the relative error in
each admissible block both in the spectral and the Frobenius norm.
The next subsection shows that we can generalize our approach to local
updates without losing the optimal complexity. 

\subsection{Algorithmic Challenges of Local Updates}

Local updates $C_{|\hat t_0\times\hat s_0} \gets C_{|\hat t_0\times\hat s_0}
+ X Y^T$ of submatrices defined by a block $b_0=(t_0,s_0)\in\ctIJ$
pose a number of additional challenges in comparison with
the global update discussed above.
In order to obtain linear complexity with respect to the size of the
local block, the top-down procedure of computing the weight matrices
and the update of coupling matrices need to be investigated more
closely.
The first one requires the weight matrix of the father and so of all
predecessors.
The second task has to update all coupling matrices even if they are not
in the sub-block of the update.

We go through four parts of the local update and discuss the special
issues:
the computation of the weight matrices, the construction of the
adaptive cluster bases for $C + X Y^T$, the update of the
${\mathcal H}^2$-matrix, and the preparation of auxiliary data
required for further updates.

The efficient computation of the weight matrix of a cluster requires
the weight matrix of the father.
If we compute an update for the root this poses no problem, but
computing the weight for a higher-level cluster would require us to
visit all of its predecessors and therefore lead to undesirable terms
in the complexity estimate.
We solve this problem by computing the weight matrices for all
clusters in a preparation step.
This can be done in linear complexity once before we start the
$LDL^T$ factorization.
For the local update we only have to recompute the weight matrices
in the sub-block of the update.
Outside of the sub-block, the matrix remains unchanged, therefore
we do not have to update the weight matrices.

There is a second challenge arising from the computation of the weight
matrices.
The blocks $(t,s_i)$ corresponding to the matrices $C_{|\hat t, \hat s_i}$
do not necessarily belong to the sub-block of the local update.
Thus we need access to all admissible blocks $(t,s_i)$ with row cluster $t$.
This is handled by lists containing all row and column blocks connected
to clusters.

As shown in subsection \ref{subsec: Adaptive Cluster Basis} the computation
of the adaptive cluster basis is a bottom-up algorithm that can be applied
to the subtree corresponding to the update.
The cluster basis outside of this subtree remains unchanged.
All predecessors can be updated by simply modifying the transfer matrix
connecting the root of the subtree to its father.
Hence there are no special problems for the local update in comparison
to the global update.

The third step is more challenging than the second one.
The coupling matrices have to be updated for all blocks $(t, s_i)$, i.e.,
they have to be multiplied by the basis change matrix $R_t$.
Since $s_i$ may lie outside of the subblock that is being updated,
we again make use of the block lists mentioned before.
In each of these blocks we only have to multiply the small matrices
$R_t$ and $S_{t, s_i}$.
Assuming again that the block tree is $c_{\rm sp}$-sparse, for one cluster
$t\in\ctI$ not more than $c_{\rm sp}$ such products have to be computed,
so the number of operations is in ${\mathcal O}(k^3)$ for one cluster.
Updating all blocks connected to the sons of $t_0$ or $s_0$ requires
${\mathcal O}(k^3 (\#\ctt+\#\cts))$ operations, where $\ctt$ and
$\cts$ denote the subtrees of $\ctI$ and $\ctJ$ with roots $t_0$
and $s_0$.
Using again the standard assumptions $\#\ctt \lesssim \#\hat t_0/k$ and
$\#\cts \lesssim \#\hat s_0/k$, we obtain a complexity of
${\mathcal O}(k^2 (\#\hat t_0 + \#\hat s_0)$.

To conclude the local update, we have to ensure that the weight matrices
are correct by recomputing them in the subtree $\ctt$ and $\cts$.
The weight matrices do not change for clusters outside the sub-block.

Altogether we end up with computational complexity in
${\mathcal O}(k^2(\#\hat t_0 + \#\hat s_0))$ for the local update in a
sub-block $b_0=(t_0,s_0)$.
Using this estimate, we can prove \cite{BORE14} that the matrix
multiplication and other higher arithmetic functions require not more
than ${\mathcal O}((\nI+\nJ) k^2 \log(n))$ operations.

%
%
\section{Slicing the Spectrum}
\label{sec:slicing}

In order to use our efficient matrix-arithmetic operations to solve an
eigenvalue problem, we use the slicing-the-spectrum algorithm that has been
described in \cite{Par80}.  For the related $\mathcal{H_{\ell}}$-matrices, which
are $\mathcal{H}$-matrices with a particularly simple block tree, the algorithm
has been investigated in \cite{BenM11c}.  Further, in \cite{BenM11c} it has been
shown by numerical examples that a generalization of the approach to
$\mathcal{H}$-matrices does not lead to an efficient algorithm in general.

We are computing the eigenvalues of a symmetric matrix. Thus all eigenvalues are
real and the function $\nu(\sigma)=\#\{\lambda\in\Lambda(A)
\ :\ \lambda<\sigma\}$ is well defined for all $\sigma\in\bbbr$.
If $\nu(a)<m\leq\nu(b)$, we know that the interval $[a,b]$ contains
the $m$-th smallest eigenvalue $\lambda_m$ of $A$.
We can run a bisection algorithm on this interval until the interval
is small enough.
The midpoint of the interval is then taken as approximation
$\hat{\lambda}_m$ of the desired eigenvalue.
We bisect the interval by computing $\nu(\tfrac{a+b}{2})$.
If $\nu(\tfrac{a+b}{2})>m$, we continue with $[a,\tfrac{a+b}{2}]$,
otherwise with $[\tfrac{a+b}{2},b]$.
We stop the algorithm if $b-a<\epsilon_{\rm ev}$ holds.
In our computations, we choose $\epsilon_{\rm ev} = 10^{-5}$.

It remains to explain how we get the inertia or $\nu(\sigma)$.
The inertia is invariant under congruence transformations, thus the
matrix $D$ of the $LDL^T$ factorization of $A$ has the same
inertia as $A$ itself.
To get $\nu(\sigma)$ we compute the $LDL^T$ factorization of
$A-\sigma I = L_{\sigma}D_{\sigma}L_{\sigma}^{T}$ and simply count the
negative diagonal entries of $D$.

\subsection{Accuracy}
\label{subsec:accuracy}

By using $\mathcal{H}^{2}$-matrices, the computation of an $LDL^T$
factorization is comparably cheap, taking essentially
${\mathcal O}(n k^2 \log n)$.
This allows the fast computation of the inertia, which would be in
${\mathcal O}(n^{3})$ for general dense matrices.
The price we have to pay is that the factorization is only approximative,
i.e., $A-\sigma I \approx \widetilde L \widetilde D \widetilde L^T$, so
we have to ensure that it is sufficiently accurate to yield the correct
value $\nu(\sigma)$.
In \cite{Par80} it is shown that this is the case if
$\|H_{\sigma}\|\leq\min_{j}|\lambda_{j}(A)-\sigma|$,
with $H_{\sigma} := (A-\sigma I) -
\widetilde{L}_{\sigma}\widetilde{D}_{\sigma}\widetilde{L}_{\sigma}^{T}$.

Thus we need a bound for the error of the form
$\|A-\widetilde{L}\widetilde{D}\widetilde{L}^{T}\| \leq \delta \|A\|$.
We further need this bound for all shifted matrices $A-\sigma I$.
In the literature the LU-decomposition has received much more attention
than the $LDL^T$ factorization.
Since the inertia of $A-\sigma I$ can also be obtained from an
LU-decomposition, we will cite some results for LU-decompositions for
hierarchical matrices:
to our best knowledge such a bound is currently not available in the
literature on $\mathcal{H}$- and $\mathcal{H}^{2}$-matrices.
In \cite{Beb05} it was shown that for certain $\mathcal{H}$-matrices
originating from certain finite element discretizations there exist
$\mathcal{H}$-mat\-rices $\widetilde{L}$ and $\widetilde{U}$ so that
$\|A-\widetilde{L}\widetilde{U}\| \leq \delta \|A\|$.
This result has been generalized in \cite{Beb07,BGK09} and more
recently in \cite{FauMP13}.
Unfortunately, it has so far not been shown that the algorithms actually
used to compute approximations yield results satisfying similar
estimates.
Fortunately, many numerical experiments show that the algorithms for
the computation of the $\mathcal{H}$-$LU$-decomposition are very good.

For the case of $A-\sigma I$, with $\sigma\neq 0$, the picture is not
positive. In \cite[Table 4.1]{BenM11c} one can see that using shifts near
eigenvalues leads to high local block ranks, which make the
$\mathcal{H}$-$LDL^T$ factorization expensive. We do not observe a similar
behavior for $\mathcal{H}^{2}$-arithmetic, but we cannot provide theoretical
bounds for the ranks.

\subsection{Generalized Eigenvalue Problem}

For the solution of generalized eigenvalue problems we have to compute the
inertia of $A-\sigma B$ instead of $A-\sigma I$.
If we think of a finite element discretization as a basis for the generalized
eigenvalue problem, then we observe that the structures of the mass and
the stiffness matrix are similar enough to allow for a cheap computation
of $A-\sigma B$ in the $\mathcal{H}^{2}$-arithmetic.
The mass matrix $B$ can be stored as a sparse matrix.
Fortunately, the nonzero entries in $B$ correspond with inadmissible
leaves in $A$, which are stored as dense matrices.
Thus the subtraction $A-\sigma B$ affects only these inadmissible leaves.

Further research should investigate the numerical properties of the
$LDL^T$ factorization of $A-\sigma B$.

\subsection{Parallelization}

The slicing of disjoint intervals is independent, thus we can easily
parallelize the algorithm by giving each node/core an instance of the
matrix and an interval to slice.
Since the size of the sparse matrix grows only linearly with the dimension of
the matrix, this is possible for comparably large matrices.
This simple parallelization has been used in \cite{BenM11c} for the
slicing algorithm for ${\mathcal H}_\ell$-matrices.
In \cite{Mac12} a speedup of 267 by using 384 processes has been reported
for a MPI-based parallelization of the algorithm from \cite{BenM11c}.
For this parallelization a master-slave structure is used.
The master provides each slave with a small interval, which the slaves
slices until all eigenvalues are found.
For these intervals the master provides a lower bound and an upper bound
and the number of eigenvalues to be found.
To provide this information some initial computations of $\nu(\sigma)$
are necessary.
These are also performed by the slaves.
The time required for the slicing of one interval varies and thus the
intervals are chosen small enough to allow for a load balancing.

This parallelization works best for many cores.
If the number of processes is small, the master process is fre\-quent\-ly
just waiting for answers, thus running 5 process on the quad-core CPU is
improving the overall run-time.

\section{Numerical Experiments}
\label{sec:numerexp}

%
%
\begin{figure}[htbp]
  \begin{center}
  \subfloat[Unit square.]{
    \includegraphics[width=0.3\textwidth]{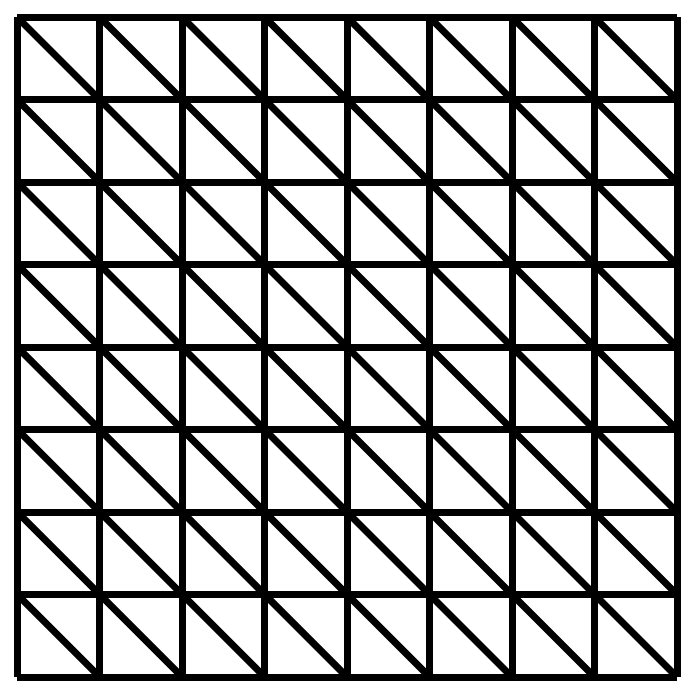}
    \label{subfig:unitsquare}
  }\subfloat[Unit circle.]{
    \includegraphics[width=0.3\textwidth]{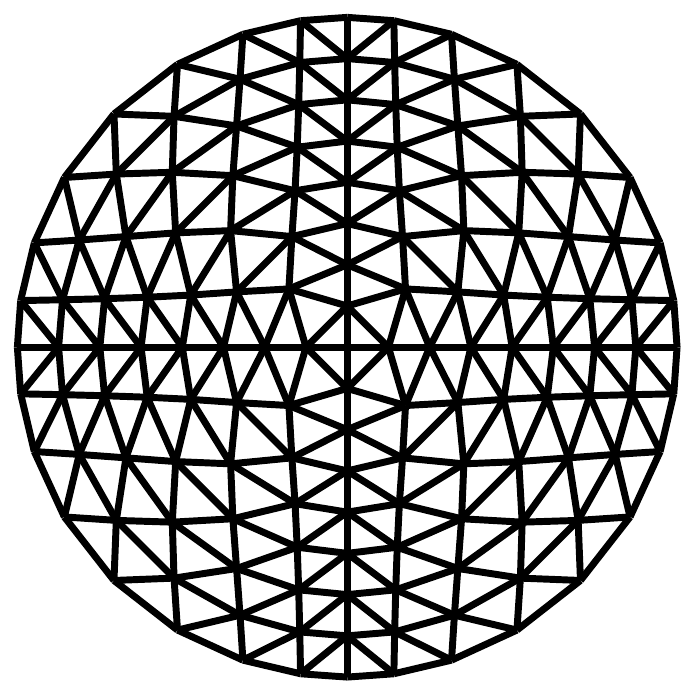}
  }

  \subfloat[L-shape.]{
    \includegraphics[width=0.3\textwidth]{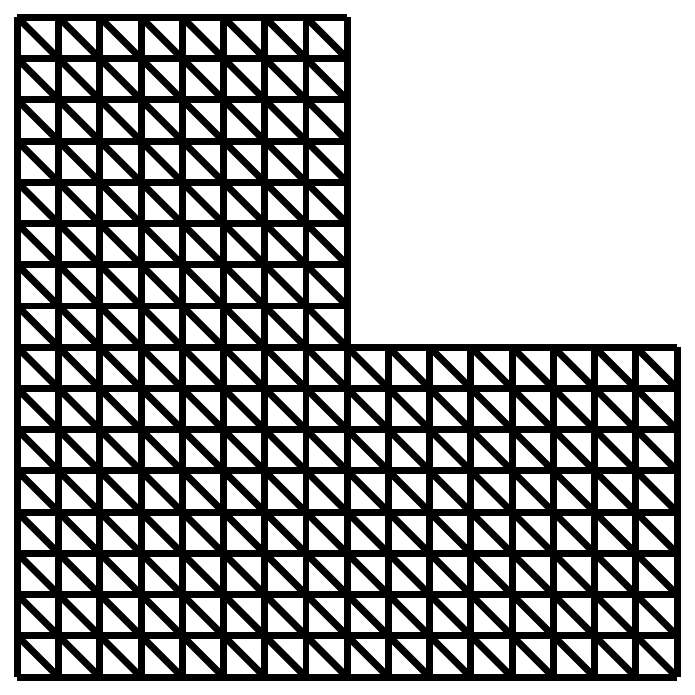}
  }\subfloat[U-shape.]{
    \includegraphics[width=0.3\textwidth]{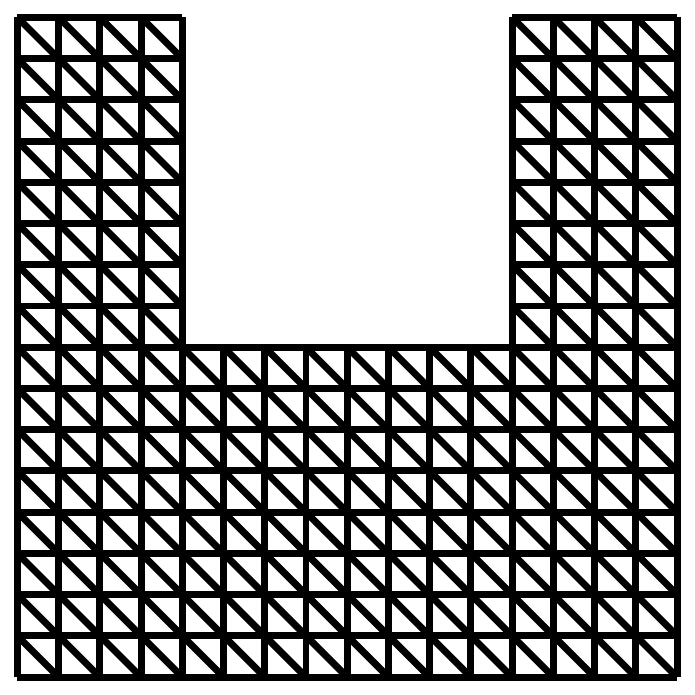}
  }

  \caption{Meshes for different geometries.} 
  \label{fig:meshes}
  \end{center}
\end{figure}

Due to the facts described in Subsection~\ref{subsec:accuracy} we cannot prove
that the proposed algorithm is accurate and efficient; at least at the moment.
Thus, numerical experiments are the only way to provide evidence that the
slicing algorithm is performing well.
For the numerical experiment we use the software
package H2Lib developed by the Scientific Computing Group at Kiel University.
This library provides examples of finite element discretizations on different
triangle meshes, see Figure~\ref{fig:meshes}.  These meshes can be refined as
needed. We use a hexa-core CPU, Intel Xeon E5645 (running at 2.40 GHz).

\begin{figure}
\begin{center}
\includegraphics[width=0.8\textwidth]{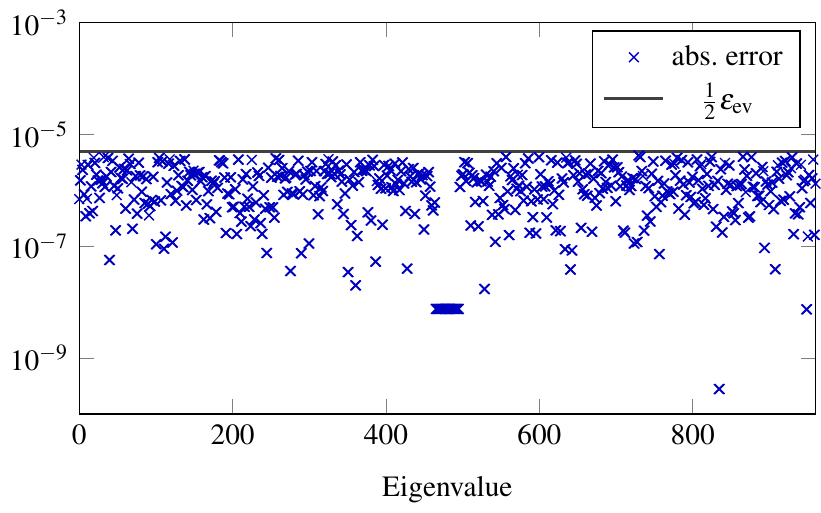}
\caption{Absolute error $|\lambda_{i}-\hat{\lambda}_{i}|$ for a
  $961\times 961$ finite element stiffness matrix corresponding to the unit square,
  $\epsilon_{\text{ev}}=10^{-5}$.}
\label{fig:abserr}
\end{center}
\end{figure}

First, we want to show that the absolute accuracy of the computed eigenvalues is
acceptable.  Therefore we use the finite element matrix related to the meshed
unit square.  We refine the mesh in
Figure~\ref{fig:meshes}\subref{subfig:unitsquare} twice, compute the eigenvalues
of this standard eigenvalue problem with the slicing algorithm and compare them
to the actual eigenvalues, which are known exactly.  In Figure~\ref{fig:abserr}
the accuracy of computed eigenvalues is shown.  The computed eigenvalues lie all
within the computed intervals.

%
%

\begin{figure}
\begin{center}
\includegraphics[width=0.8\textwidth]{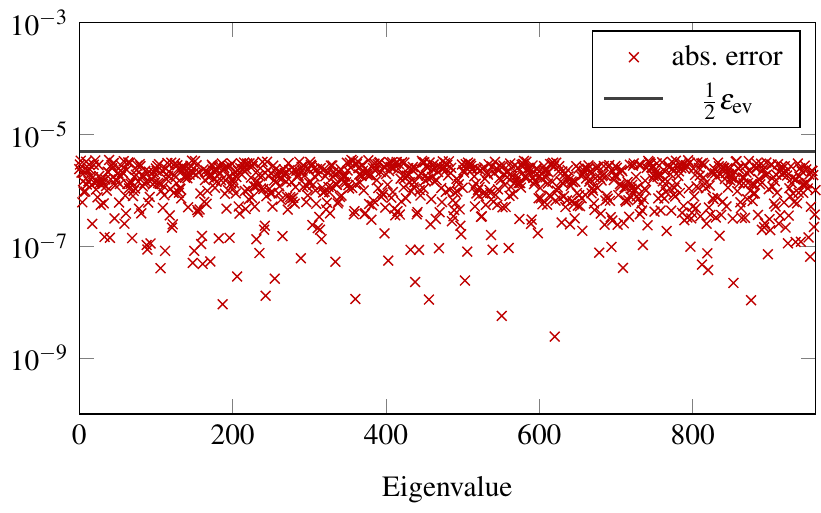}
\caption{Absolute error $|\lambda_i-\tilde\lambda_i|$ for the
  generalized eigenproblem for a $961\times 961$ finite element
  matrix corresponding to the unit square, $\epsilon_{\text{ev}}=10^{-5}$.}
\label{fig:abserrgen}
\end{center}
\end{figure}

On the same mesh we then compute the mass matrix and solve the generalized
eigenvalue problem, both with the \lapack~\cite{Lapack} eigenvalue solver for
symmetric generalized eigenvalue problems \texttt{dsygv} and with the
slicing algorithm.
The result is similar to the previous one, as we observe in
Figure~\ref{fig:abserrgen} that again the allowed tolerance is fulfilled
for all eigenvalues.

%
%

\begin{figure}\begin{center}
\includegraphics[width=0.8\textwidth]{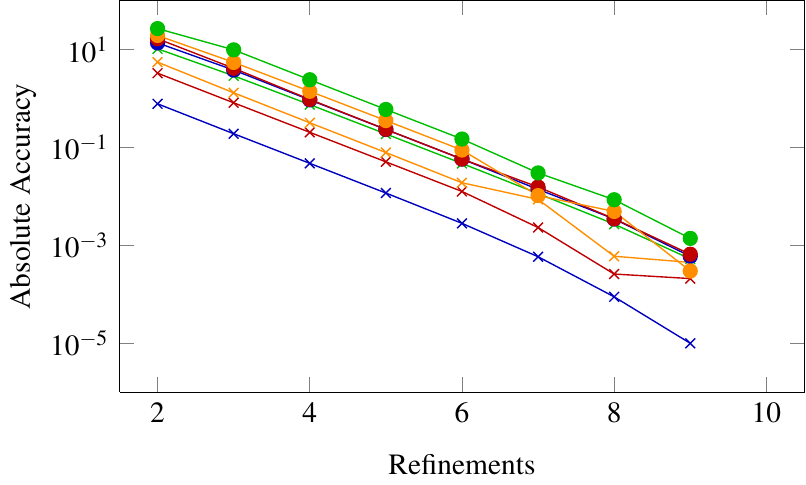}
\caption{Convergence of the eigenvalues with respect to the mesh parameter.}
\label{fig:convergence}
\end{center}
\end{figure}

Since we are solving finite element eigenvalue problems, we expect the
smallest eigenvalues to converge to the eigenvalues of the differential
operator.
This can be seen in Figure~\ref{fig:convergence} for the 8 smallest
eigenvalues, where 3 refinements correspond to the mesh shown in
Figure~\ref{fig:meshes}\subref{subfig:unitsquare}:
we obtain the ${\mathcal O}(h^2)$ convergence predicted by standard
theory.

%
%

\begin{figure}
\begin{center}
\includegraphics[width=0.8\textwidth]{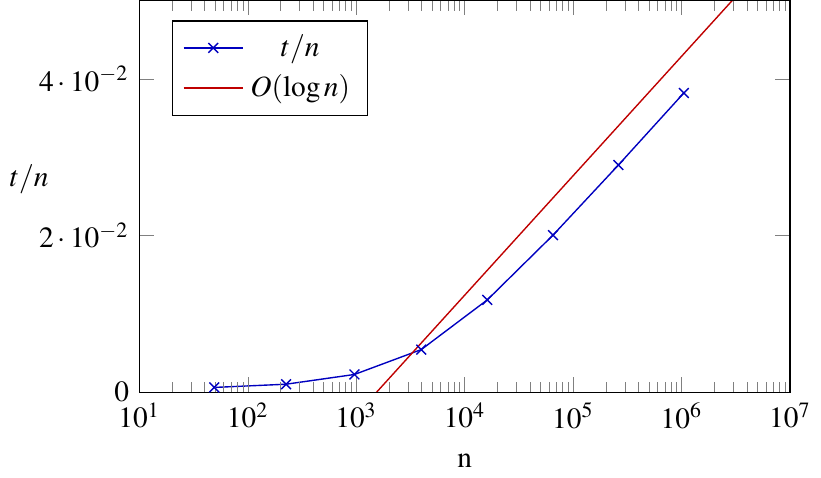}
\caption{Runtime divided by matrix dimension; unit circle, 8 smallest
  eigenvalues of the generalized eigenvalue problem.}
\label{fig:runtime}
\end{center}
\end{figure}

In Table~\ref{tab:runtime_eight} the runtime, the time for one slice, and the
accuracy are shown for different refinements of the meshes in
Figure~\ref{fig:meshes}.  The accuracy is the maximum absolute error for the
computed eigenvalues compared with the results from the \lapack{} eigensolver
\texttt{dsygv}.  For matrices with $n\geq 5000$ the accuracy is not computed,
since the dense matrices are too large and the computations with \lapack{} would
take too long.  Figure~\ref{fig:runtime} shows the time per degree of freedom
using a logarithmic scale for $n$.  It seems to suggest a complexity of
${\mathcal O}(n \log n)$ for large values of $n$, i.e., the effective rank $k$
of the ${\mathcal H}^2$-matrix approximation of the $LDL^T$ factorization
appears to be bounded independently of the mesh size.

%
%
\begin{table}
  \centering
  \begin{tabular*}{0.8\textwidth}{@{\extracolsep{\fill}}rrrr@{}}
    \toprule
    \multicolumn{4}{@{}c}{Unit Square}\\
    $n$ & $t_{\text{8 ev}}$ in s & $t_{\text{single slice}}$ in s & maximal abs. err. \\
    \midrule
          $225$ &      $0.23$ &        $<0.01$ &        $3.1510_{-06}$ \\
          $961$ &      $2.18$ &        $0.01\pm 0.00$ & $2.8435_{-06}$ \\
        $3,969$ &     $21.61$ &        $0.14\pm 0.01$ & $3.2226_{-06}$ \\
       $16,129$ &    $190.20$ &        $1.30\pm 0.02$ & \\
       $65,025$ &  $1,304.22$ &        $9.31\pm 0.10$ & \\
      $261,121$ &  $7,577.24$ &       $56.08\pm 0.62$ & \\
    $1,046,529$ & $39,992.06$ &      $305.03\pm 5.18$ & \\
    \midrule
    \multicolumn{4}{@{}c}{Unit Circle}\\
    $n$ & $t_{\text{8 ev}}$ in s & $t_{\text{single slice}}$ in s & maximal abs. err. \\
    \midrule
        $481$ &      $0.79$ &   $0.01\pm 0.01$ & $2.3895_{-06}$ \\
      $1,985$ &      $9.02$ &   $0.06\pm 0.01$ & $2.6512_{-06}$ \\
      $8,065$ &     $93.70$ &   $0.68\pm 0.02$ & \\
     $32,513$ &    $829.55$ &   $5.92\pm 0.11$ & \\
    $130,561$ &  $5,457.30$ &  $41.30\pm 0.59$ & \\
    $523,265$ & $32,327.59$ & $252.34\pm 2.75$ & \\
    \midrule
    \multicolumn{4}{@{}c}{L-Shape}\\
    $n$ & $t_{\text{8 ev}}$ in s & $t_{\text{single slice}}$ in s & maximal abs. err. \\
    \midrule
    $161$ &      $0.14$ &          $<0.01$ & $1.8366_{-06}$ \\
    $705$ &      $1.42$ &   $0.01\pm 0.00$ & $3.8873_{-06}$ \\
    $2,945$ &     $12.33$ &   $0.07\pm 0.01$ & $3.8643_{-06}$ \\
    $12,033$ &    $125.62$ &   $0.69\pm 0.01$ & \\
    $48,641$ &    $999.53$ &   $5.52\pm 0.06$ & \\
    $195,585$ &  $6,587.22$ &  $36.37\pm 0.37$ & \\
    $784,385$ & $27,722.95$ & $152.72\pm 5.84$ & \\
    \midrule
    \multicolumn{4}{@{}c}{U-Shape}\\
    $n$ & $t_{\text{8 ev}}$ in s & $t_{\text{single slice}}$ in s & maximal abs. err. \\
    \midrule
    $153$ &      $0.14$ &          $<0.01$ & $2.0762_{-06}$ \\
    $689$ &      $1.82$ &   $0.01\pm 0.00$ & $3.7008_{-06}$ \\
    $2,913$ &     $10.14$ &   $0.06\pm 0.01$ & $3.0300_{-06}$ \\
    $11,969$ &     $99.89$ &   $0.56\pm 0.01$ & \\
    $48,513$ &    $808.08$ &   $4.61\pm 0.04$ & \\
    $195,329$ &  $5,607.98$ &  $31.99\pm 0.36$ & \\
    $783,873$ & $24,464.86$ & $139.24\pm 5.04$ & \\
    \bottomrule
  \end{tabular*}
  \caption{Runtime for the computation of the 8 smallest eigenvalues on
    different shapes; for small matrices including the accuracy.}
  \label{tab:runtime_eight}
\end{table}

In Table~\ref{tab:comparewith4andlapack} we compare the algorithm with the
slicing algorithm for $\mathcal{H}$-matrices described in \cite{BenM11c}. Since
the $\mathcal{H}$lib \cite{hlib} is more optimized with respect to speed than
the H2lib we choose to reimplement the algorithm from \cite{BenM11c} in the
H2lib for a fair comparison.
Thereby we also generalized the algorithm to generalized eigenvalue
problems. We see that the implementation based on $\mathcal{H}^{2}$ is slightly
faster at the same accuracy.

However, using \lapack{} \texttt{dsygv}, based on an implicit QZ algorithm on
the dense matrix, would be much faster for the computation of all eigenvalues. A
backward stable algorithm is used to compute the eigenvalues to almost machine
precision. The generalized eigenvlaue problem, unit-square with mass matrix, of
dimension 3969 can be solved in 38 s and the problem of dimension 16129 in 2296
s. The bigger problem requires about 2 GB storage. Thus one should only use the
slicing algorithm for large problems.

%
%
\begin{table}
  \centering
  \begin{tabular*}{0.8\textwidth}{@{\extracolsep{\fill}}rrrrrr@{}}
    \toprule
    \multicolumn{6}{@{}c}{Unit Square, with mass matrix}\\
    $n$ & ev 
    & $t_{\mathcal{H}}$ &
    error & $t_{\mathcal{H}^2}$ & error \\
    \midrule
       225 & 8 
           &   0.24 & $3.15_{-06}$%
           &   0.23 & $3.15_{-06}$ \\%
       961 & 8 
           &   2.60 & $2.84_{-06}$%
           &   2.18 & $2.84_{-06}$ \\%
     3,969 & 8 
           &  29.07 & $3.22_{-06}$%
           &  21.61 & $3.22_{-06}$ \\%
    16,129 & 8 
           & 246.83 & ---%
           & 190.20 & --- \\
    65,025 & 8 
           & 1,555.23 & ---%
           & 1,304.22 & ---\\%
    261,121 & 8 
            &          & ---%
            & 7,577.24 & ---\\%

    \bottomrule
  \end{tabular*}
  \caption{Comparing the $\mathcal{H}^{2}$ slicing algorithm with
    the $\mathcal{H}$ slicing algorithm described in \cite{BenM11c}. All timings in
    seconds.}
  \label{tab:comparewith4andlapack}
\end{table}

%
%
\begin{table}
  \centering
  \begin{tabular*}{0.8\textwidth}{@{\extracolsep{\fill}}lccr@{}}
    \toprule
    \multicolumn{4}{@{}c}{Unit Square}\\
    $n$ & no. of cores & $t_{\text{all ev}}$ in s & speedup \\
    \midrule
    961 & 1+0 & 185.40 & single core code\\
    961 & 1+1 & 182.73 & 1.01\\
    961 & 2+1 & \phantom{1}99.85 & 1.86\\
    961 & 3+1 & \phantom{1}68.71 & 2.70\\
    961 & 4+1 & \phantom{1}47.77 & 3.88\\
    961 & 5+1 & \phantom{1}38.60 & 4.80\\
    \bottomrule
  \end{tabular*}
  \caption{Speedup by parallelization; generalized eigenvalue problem, all eigenvalues.}
  \label{tab:parallel}
\end{table}

%
%

Finally we test the MPI based parallelization, see Table~\ref{tab:parallel}.
Here we use a quadcore CPU, Intel Core i5-3570 (running at 3.40 GHz) and compute
the speedup in comparison with the runtime of the single core code.  Since the
master is not doing any work we see good speedups for up to 4 slave processes.

\section{Conclusions}

We have investigated whether the computation of eigenvalues of symmetric
$\mathcal{H}^{2}$-matrices can be done efficiently by slicing the spectrum.  Our
results show that for small $n$ other methods, eventually even dense eigenvalue
solver, are more efficient.  However, the experiments further show that the
computational costs per eigenvalue scale with $\mathcal{O}(n \log n)$ and thus
for large $n$ the method will be very efficient. It remains open whether the
usage of $\mathcal{H}^{2}$-arithmetic is significantly more efficient than
$\mathcal{H}$-arithmetic or not. The additional structure might be used for
higher efficiency, but produces also more overhead. 

\bibliographystyle{plain}
\bibliography{bib}

\end{document}